\documentclass[a4paper,twoside]{article}
\usepackage{amsmath,amsfonts,amssymb,amsthm}
\usepackage{enumitem}
\usepackage[bookmarks]{hyperref}
\usepackage[center]{titlesec}
\usepackage{
mathrsfs,dsfont}
\usepackage{tocbibind}
\usepackage{geometry}
\numberwithin{equation}{section}

\newtheorem{thm}{Theorem}[section]
\newtheorem{cj}{Conjecture}[section]
\newtheorem{defi}[thm]{Definition}
\newtheorem{lem}[thm]{Lemma}
\newtheorem{cor}[thm]{Corollary}
\newtheorem{rem}{Remark}

\def\bthm{\begin{thm}\def\ethm{\end{thm}}}

\def\blem{\begin{lem}\def\elem{\end{lem}}}

\def\bitm{\begin{itemize}} \def\eitm{\end{itemize}}
\def\benu{\begin{enumerate}} \def\eenu{\end{enumerate}}

\def\bpf{\begin{proof}}\def\epf{\end{proof}}
\def\beq{\begin{equation}}\def\eeq{\end{equation}}
\def\beqs{\begin{eqnarray}}\def\eeqs{\end{eqnarray}}
\def\beqsnl{\begin{eqnarray*}}\def\eeqsnl{\end{eqnarray*}}

\def\bb{\mathbb}
\def\ca{\mathcal}

\def\scr{\mathscr}

\newcommand{\bra}[1]{\left(#1\right)}

\newcommand{\abs}[1]{\left\vert#1\right\vert}
\newcommand{\set}[1]{\left\{#1\right\}}

\newcommand{\norm}[1]{\left\Vert#1\right\Vert}

\DeclareMathOperator{\re}{Re}
\DeclareMathOperator{\im}{Im}

\usepackage{fancyhdr}
\begin{document}
\title{\textbf{
Sharp 
Hardy-Littlewood-Sobolev Inequalities on Octonionic Heisenberg Group}}
\footnotetext{\emph{Date:} Feb. 20th, 2014, revised on July 12, 2014.}
\footnotetext{Michael Christ is supported in part by NSF grant DMS-0901569. Heping Liu is supported by National Natural Science Foundation of China under Grant \#11371036 and
the Specialized Research Fund for the Doctoral Program of Higher Education of China under Grant \#2012000110059. An Zhang is supported by China Scholarship Council under Grant No. 201306010009.}

\author{Michael Christ, Heping Liu, An Zhang\footnote{Corresponding author.}}

\date{}
\maketitle
\begin{abstract}
This paper is a second one following \cite{clz} in series, considering sharp Hardy-Littlewood-Sobolev inequalities on groups of Heisenberg type. The first important breakthrough was made in \cite{frank2012sharp}. In this paper, analogous results are obtained for octonionic Heisenberg group.
\end{abstract}

\section{Introduction}
The Hardy-Littlewood-Sobolev inequalities (HLS) for conjugate exponent on group of Heisenberg type states that
\[|\iint_{G\times G}\frac{\overline{f(u)}g(v)}{|u^{-1}v|^\lambda}dudv| \lesssim  \|f\|_p \|g\|_p,\] with $0<\lambda<Q, p=\frac{2Q}{2Q-\lambda}, Q$ is the homogeneous dimension of the group $G$.  In \cite{frank2012sharp}, they proved on Heisenberg group that for all $\lambda$, extremizer  for above HLS inequality is almost uniquely \[\left((1+|z|^2)^2+|t|^2\right)^{-\frac{2Q-\lambda}{4}}\] with group elements parameterized by $u=(z,t), z\in \bb{C}^n, t\in \bb{R}$. In \cite{clz}, we extended this result to quaternionic Heisenberg group, which states that for partial exponent $\lambda\ge 4$, extremizer is almost uniquely of the similar form.  This paper duplicates the method in above two works and we proved on the only 15-dimensional octonionic Heisenberg group that for $\lambda\ge 12$, extremizer for HLS inequality exists almost uniquely of the form \[\left((1+|z|^2)^2+|t|^2\right)^{\frac{\lambda}{4}-11}\] with group elements $u=(z,t), z\in \bb{O}, t\in \im \bb{O}$, where $\bb{O}$ is the octonions. Note that the homogeneous dimension of this octonionic Heisenberg group is 22, so this result coincides with the cases on above classical Heisenberg and quaternionic Heisenberg groups. Dual conformally-invariant Sobolev inequalities (associated with intertwining operators of complementary representation of exceptional Lie group $F_{4(-20)}$) and endpoint limit case, Log-Sobolev inequality, are also obtained. As in \cite{frank2012sharp,clz}, we still use Cayley transform to move the inequalities on group onto sphere, then it's natural to see that \emph{zero center-mass} technique is a perfect way to break the huge conformal symmetry group of the inequalities, which is a critical property of HLS-type inequalities.

\section{Sharp HLS and Related Sobolev-type Inequalities}
The only octonionic Heisenberg group is a 15-dimensional 2-step nilponent group parameterized by its Lie algebra $G=\bb{O}\times \im \bb{O}$ with group law
\[u  u'=(z,t) (z',t')=(z+z',t+t'+2\im z \overline{z'}),\] where $u=(z,t)$ is the group element with $z\in \bb{O}, t\in \im \bb{O}$. Here $\bb{O}$ is octonions, the noncommutative (also non-associative) 8-dimensional division ring over real numbers, $\overline{z'}$ is the conjugate and $z\overline{z'}$ is the standard multiplication in $\bb{O}$. This group can be viewed as the nilponent part of Iwasawa
decomposition of rank one connected Lie group $F_{4(-20)}$, the exceptional Lie group, which is the isometry group of rank one hyperbolic symmetric
space over octionions, which can also be identified with homogeneous space $F_{4(-20)}/Spin(9)$. The symmertic space can be realized in two models: Siegel domain model $D$ and unit ball model $B$ with boundary respectively identified with $G$ and octonionic sphere $S\simeq\bb{S}^{15}$. Cayley transform will give a bijection between the two models. ``Boundary" Cayley transform is defined by
\begin{align}
\ca{C}:~~\qquad G ~&\longrightarrow~ S\setminus \{o\}\nonumber\\
u=(z,t) ~&\longmapsto~ \zeta=(\zeta_1,\zeta_2)=\left(\frac{2z}{1+|z|^2-t},\frac{1-|z|^2+t}{1+|z|^2-t}\right)\nonumber\\
\ca{C}^{-1}: ~S\setminus \{o\} ~&\longrightarrow~ G\nonumber\\
\zeta=(\zeta_1,\zeta_2) ~&\longmapsto~ u=(z,t)=\left(\frac{\zeta_1}{1+\zeta_2}, -\text{Im}\frac{1-\zeta_2}{1+\zeta_2}\right)
\end{align} with $z,\zeta_1,\zeta_2\in \bb{O}, t\in \im \bb{O}, o=(0,-1)$ the south pole on $S$ and Jacobian determinant
\begin{align*}
|J_{\ca{C}}(u)|= & 2^{Q-7}\big((1+|z|^2)^2+|t|^2\big)^{-\frac{Q}{2}}\\
= & 2^{-7}|1+\zeta_{n+1}|^{Q}.
\end{align*}
This ``boundary" Cayley transform is a generalization of stereographic projection in Euclidean space.
We denote the homogenous norm on $G$ \[ |u|=|(z,t)|=(|z|^4+|t|^2)^{\frac{1}{4}}, \] and the group distance of
two elements is then defined by
\[d_G(u,v)=|v^{-1}u|= \left(|z-z'|^4+|t-t'+2\im z\overline{z'}|^2\right)^{\frac{1}{4}}.\]
We use notation $Q=8+7\times2=22$ for the homogeneous dimension of group $G$, $\delta u=(\delta z, \delta^2t)$ for group dilation
and $du=dzdt$ for Haar measure with $dzdt$ the Lebesgue measure on $\bb{R}^8\times \bb{R}^7$. We denote $Aut(G)$ the set of all conformal
transformations on $G$ (diffeomorphisms preserving the contact structure, also called octonionic automorphisms, that's why we use the notation), element of which is composition of translations, rotations (on $z$),
dilations and inversion
\beq\label{inversion}\sigma_{inv}:(z,t)\mapsto (-\frac{z}{|z|^2-t},-\frac{t}{|z|^4+|t|^2}),\eeq
extension of which on related rank one symmetric space $F_{4(-20)}/Spin(9)$ is an isometry.

In this paper, we use notation $``\sim"$ for $``="$, ignoring a constant multiple independent of functions or what we care in formulas. Sharp HLS inequality on $G$ is then given in the following theorem.

\bthm\label{t-hls} \emph{[{Sharp HLS on Octonionic Heisenberg Group}]}\\
Let $12\le \lambda<Q=22, p=\frac{2Q}{2Q-\lambda}$, then $\forall ~f,g\in L^p(G)$,
\beq\label{hls}
\Big|\iint_{G\times G}\frac{\overline{f(u)}g(v)}{|u^{-1}v|^\lambda}dudv\Big| \le C_\lambda \|f\|_p\|g\|_p
\eeq with sharp constant
\beq
C_\lambda
=~2^{-\frac{4\lambda}{Q}}\bra{\frac{2\pi^8}{7!}}^{\frac{\lambda}{Q}}
\frac{7!\Gamma(\frac{Q-\lambda}{2})}{\Gamma(\frac{2Q-\lambda}{4})\Gamma(\frac{2Q-\lambda}{4}-3)}.
\eeq
Moreover, all extremizers are given by
\beq\label{e} f\sim g \sim (|J_{\ca{C}}\circ\sigma||J_\sigma|)^{\frac{1}{p}}
\sim\left||z|^2+t-2z_0\overline{z}+w_0\right|^{-\frac{2Q-\lambda}{2}},\eeq
with $\sigma\in Aut(G),
z_0, w_0\in \bb{O}$, satisfying $\emph{Re}w_0>|z_0|^2$, and we can choose $\sigma=\ca{S}_{\delta_0}\circ\ca{L}_{u_0}$ or $\sigma=\ca{S}_{\delta}\circ\ca{C}^{-1}\circ A_\xi\circ \ca{C}$, where $\ca{S}_{\delta_0}, \ca{S}_\delta$ are group dilations, $\ca{L}_{u_0}$ is group left translation and
$A_\xi$ is a sphere rotation in $Spin(9)$ s.t. $A_\xi^{-1}(0,1)=\frac{\xi}{|\xi|}$ with $(0,1)$ the north pole on $S$, with parameters $\delta_0=(\re w_0-|z_0|^2)^{-\frac{1}{2}}, u_0=(z_0,-\im w_0)$, $\xi=(\frac{2z_0}{w_0+1},\frac{w_0-1}{w_0+1}), \delta=\sqrt{\frac{1\mp|\xi|}{1\pm|\xi|}}$.
\ethm

Using the Cayley transform, we can give the equivalent sphere editon of last theorem.
If we define the distance function on $S$ to be
\[d_S(\zeta,\eta)=2^{-\frac {1}{2}}|1-\zeta\cdot\overline{\eta}|^{\frac{1}{2}}\] with $\zeta\cdot\overline{\eta}=\zeta_1\overline{\eta_1}+\zeta_2\overline{\eta_2}$, then there exists the following relation between two distances on $G$ and $S$,
\begin{align}\label{dr}
d_S(\zeta,\eta)=& ((1+|z|^2)^2+|t|^2)^{-\frac{1}{4}}((1+|z'|^2)^2+|t'|^2)^{-\frac{1}{4}}d_G(u,v)\nonumber\\
=& 2^{\frac{7}{Q}-1}|J_{\ca{C}}(u)|^{\frac{1}{2Q}}|J_{\ca{C}}(v)|^{\frac{1}{2Q}}d_G(u,v).
\end{align}
From this relation, we obtain sharp HLS inequality on octonionic sphere,
through the correspondence between functions $f$ on $G$ and $\tilde{f}$ on $S$,
\beq\label{cor}\tilde{f}(\zeta)=f(\ca{C}^{-1}\zeta)|J_{\ca{C}^{-1}}|^{\frac{1}{p}}.\eeq
We denote $Aut(S)=\set{\tau=\ca{C}\circ \sigma \circ \ca{C}^{-1}: \sigma\in Aut(G)}$ the set of all conformal transformations on $S$, then sharp HLS inequality on octonionic sphere is given in the following theorem.
\bthm\label{t-shls} \emph{[{Sharp HLS on Octonionic Sphere}]}\\
Let $12\le \lambda<Q=22, p=\frac{2Q}{2Q-\lambda}$, then $\forall ~f,g\in L^p(S)$,
\beq\label{shls}
\Big|\iint_{S\times S}\frac{\overline{f(\zeta)}g(\eta)}{d_S^\lambda(\zeta,\eta)}d\zeta d\eta\Big| \le C'_\lambda \|f\|_p\|g\|_p
\eeq with sharp constant
\begin{align}\label{C'}
C'_\lambda=&~ 2^\frac{15\lambda}{Q}C_\lambda\nonumber\\
=&~ 2^{\frac{\lambda}{2}}\bra{\frac{2\pi^8}{7!}}^{\frac{\lambda}{Q}}
\frac{7!\Gamma(\frac{Q-\lambda}{2})}{\Gamma(\frac{2Q-\lambda}{4})\Gamma(\frac{2Q-\lambda}{4}-3)}.
\end{align}  Moreover, all extremizers are given by
\beq\label{e-s} f\sim g\sim |J_\tau|^{\frac{1}{p}}\sim|1-\xi\cdot\bar{\zeta}|^{-\frac{2Q-\lambda}{2}},\eeq
with $\tau\in Aut(S), \xi\in \bb{O}^2, |\xi|<1$, and we can choose  $\tau=\ca{C}\circ\ca{S}_{\delta_0}\circ\ca{L}_{u_0}\circ\ca{C}^{-1}$ or ~$\ca{C}\circ\ca{S}_\delta\circ\ca{C}^{-1}\circ A_\xi$,
where $\ca{S}_{\delta_0}, \ca{S}_\delta$ are group dilations, $\ca{L}_{u_0}$ is group left translation, with parameters $\delta_0=\frac{|1+\xi_2|}{\sqrt{1-|\xi|^2}}, u_0=(\frac{\xi_1}{1+\xi_2},-\im \frac{1-\xi_2}{1+\xi_2}), \delta=\sqrt{\frac{1\pm|\xi|}{1\mp|\xi|}}$, and $A_\xi$ is a rotation in $Spin(9)$ s.t.
$A_\xi^{-1}(0,1)=\frac{\xi}{|\xi|}$ with $(0,1)$ the north pole on $S$.
\ethm
We leave the proofs of sharp HLS in next section and first give several small remarks about the two theorems.\\
\textbf{Remark:}
\bitm
\item
Existence of extremizer holds for all $0<\lambda<Q$ for Theorem \ref{t-hls} and \ref{t-shls}.
Several a bit standard methods can be used to prove the existence while compactness is the basic idea. For range $\lambda<12$, unfortunately, we get no global result, but a weaker statement are proved that the optimal functions in above theorems for corresponding power $\lambda$ are \emph{local} extremizers. We omit the proof, see \cite{frank2012sharp,clz}.
\item
The large conformal symmetry group of the HLS inequality (\ref{hls}) consists of constant multiple, left-translation, dilation,
rotation (on $q$ variable) and inversion: $f(u)\mapsto f(\sigma_{inv} u) |u|^{-\frac{2Q}{p}}$, where $\sigma_{inv}$ is the
group inversion(\ref{inversion}). In other words, the inequality is invariant under the conformal action
$f\mapsto f\circ\sigma |J_\sigma|^{\frac{1}{p}}, \forall \sigma\in Aut(G)$. Similarly, HLS inequality (\ref{shls}) is invariant under the conformal action $f\mapsto f\circ\tau |J_\tau|^{\frac{1}{p}}, \forall \tau\in Aut(S)$.
\item
Theorem \ref{t-hls} tells that, modulo the conformal symmetry group of (\ref{hls}) (or partially constant multiple, left-translation and dilation),
extremizer exists uniquely,  i.e. sharp equality (\ref{hls}) holds if and only if
\[ f=g=\big((1+|z|^2)^2+|t|^2\big)^{-\frac{2Q-\lambda}{4}}.\]
\item
Theorem \ref{t-shls} tells that, modulo the conformal symmetry group of (\ref{shls})(or partially constant multiple, $\ca{C}\circ\ca{S}_\delta\circ\ca{C}^{-1}, \ca{C}\circ\ca{L}_u\circ\ca{C}^{-1}$, where $\ca{S}_\delta,\ca{L}_u$ are respectively dilation and left translation on group $G$), extremizer exists uniquely, i.e. sharp equality (\ref{shls}) holds if and only if
\[ f=g=1.\]
\eitm

Now we want to give two related sharp Sobolev-type inequalities.

We can write above sharp HLS inequalities in a dual form concerning intertwining operators, which we may call sharp
conformally-invariant Sobolev inequalities. First, we introduce and say something about intertwining operators associated to spherical principle series representation of execptional Lie group $F_{4(-20)}$. We have $Spin(9)-$irreducible decomposition of $L^2(S)$
\beq\label{dec} L^2(S)=\bigoplus_{j\ge k\ge 0} W_{j,k},\eeq where $W_{j,k}$ is ``$(j,k)-$bispherical harmonic subspaces", which is a finite dimensional space spanned by elements from the cyclic action of $Spin(9)$ on zonal harmonics
\begin{align}\label{zh}
Z_{j,k}(\zeta)\sim& \cos^{j-k}\phi~{}_2F_1\bra{-\frac{j-k}{2},\frac{1-(j-k)}{2};\frac{7}{2};-\tan^2\phi}\nonumber\\
&\times\cos^{j+k}\theta~{}_2F_1\bra{-k,\frac{-j-6}{2};4;-\tan^2\theta},
\end{align}
with $\abs{\zeta_2}=\cos\theta, \re\zeta_2=\cos\theta\cos\phi,~{}_2F_1(a,b;c;z)$ is the hypergeometric function.
For $d\in(0,Q)$, we define intertwining operator $\ca{A}_d$ on $S$ to be an operator diagonal with regard to above bispherical harmonic decomposition (\ref{dec}) with spectrums being
\beq
\ca{A}_d|_{W_{j,k}}= \frac{\Gamma(j+\frac{Q+d}{4})}{\Gamma(j+\frac{Q-d}{4})}\frac{\Gamma(k+\frac{Q+d}{4}-3)}{\Gamma(k+\frac{Q-d}{4}-3)}.
\eeq
The definition is, modulo a constant multiple, equivalent to the intertwining relationship
\[|J_\tau|^{\frac{Q+d}{2Q}}(\ca{A}_df)\circ\tau=\ca{A}_d\left(|J_\tau|^{\frac{Q-d}{2Q}}(f\circ\tau)\right), \forall \tau\in Aut(S), f\in \scr{D}(S),\] where $|J_\tau|$ is the Jacobian determinant of $\tau$, $\scr{D}(S)$ is the space of smooth function of compact support on $S$. The operators can then be extended from $\scr{D}(S)$ to Folland-Stein-Sobolev space $W^{d/2,2}(S)$. A basic result states that fundamental solution of the intertwining operator is constant multiple of power of distance function, which is given accurately by
\beq \ca{A}_d^{-1}(\zeta,\eta)=c_d d_S^{d-Q}(\zeta,\eta)\eeq
with
\beq\label{c}
c_d^{-1}=\frac{2^{\frac{Q-d}{2}+1}\pi^{\frac{Q}{2}-3}\Gamma(\frac{d}{2})}{\Gamma(\frac{Q-d}{4})\Gamma(\frac{Q-d}{4}-3)}.
\eeq
Above results were more or less proved in \cite{johnson1976composition} and (1) in Lemma \ref{l-eig}, where we compute eigenvalues of integral operators with kernel of power of distance. See \cite{branson2013moser} and its references for intertwining operators of principle series representation of semisimple Lie group and for a different method of similar proof by analysis language for Heisenberg group.
Now, we give the sharp conformally-invariant Sobolev inequality on $S$ in the following theorem. Note that similar definition and results exist for intertwining operators defined on group $G$.
\bthm \emph{[Sharp Conformally-Invariant Sobolev Inequality]}\\
Let $0<d\le Q-12, q=\frac{2Q}{Q-d}$, then $\forall f \in W^{d/2,2}(S)$,
\beq
\int_S \bar{f}\ca{A}_df\ge C_d'' \|f\|_q^2,
\eeq with sharp constant
\[C_d''=(c_dC'_{Q-d})^{-1},\]
 and all extremizers
\[f\sim |1-\xi\cdot\bar{\zeta}|^{-\frac{Q-d}{2}},\] with $\xi\in \bb{O}^2, |\xi|<1$. Constants $c_d, C'_{Q-d}$ are respectively given by \emph{(\ref{c})(\ref{C'})}.
\ethm

Also, as in \cite{frank2012sharp,clz}, we can give the endpoint limit analogue of sharp HLS inequality at $\lambda=Q$, using standard functional limit argument. The endpoint case corresponds to Log-Sobolev inequality. We list the result on sphere in the following theorem.

\bthm \emph{[Sharp Log-Sobolev Inequality]}\\
$\forall f\ge 0\in L^2\emph{Log}L(S)$, normalized by $\int_S f^2=|S|$,
\beq\label{log}\iint_{S\times S}\frac{|f(\zeta)-f(\eta)|^2}{d^Q_S(\zeta,\eta)}d\zeta d\eta \ge C \int_S f^2\log f^2
\eeq
with sharp constant \[C=
\frac{2^{\frac{Q}{2}+3}\pi^8}
{Q\Gamma(\frac{Q}{4})\Gamma(\frac{Q}{4}-3)}
,\] and some extremizers
\[ f\sim |1-\xi\cdot\bar{\zeta}|^{-\frac{Q}{2}}\] with nonzero $\xi\in\bb{O}^2, |\xi|<1$ and $f$ satisfying normalized condition.
\ethm
The proofs of conformally-invariant Sobolev and Log-Sobolev inequalities are absolutely the same as that in \cite{clz}.

\section{Proof of Sharp HLS}
\emph{Step 1:}\\
For simplicity, we prove sharp HLS on sphere. The group case can then be implied from sphere case by relation (\ref{dr}) and (\ref{cor}). Existence was proved by similar arguments in \cite{frank2012sharp,clz}.
We know that the integral kernal $d_S(\zeta,\eta)^{-\lambda}$ is positive definite.
Therefore we can restrict the sharp problem to $f=g$ case.
By standard argument, we can further restrict the extremizer to be a complex multiple of positive real-value function.
Acutally, for any $f=a+ib$, with $a,b$ respectively the real and image part function on $S$.
Then the left side of the HLS inequality (\ref{shls}) $I(f)=I(a)+I(b)$ because of the symmetry of the distance.
By Cauchy-Schwartz inequality, $I(f)\le I(|f|=\sqrt{a^2+b^2})$ with equality holds if and only if
$a(\zeta)a(\eta)\equiv b(\zeta)b(\eta)$ and $a(\zeta)a(\eta)\ge 0$, which gives that $f=(\frac{a}{b}+i)b=c|b|, c\in \bb{C}$.
So we can assume the extremizer $h$ is nonnegtive. The vanishing first variation give Euler-Lagrange equation for $h$,
\beq h^{p-1}(\zeta)\sim \int_S \frac{h(\eta)}{|1-\zeta\cdot\bar{\eta}|^{\frac{\lambda}{2}}}d\eta,\eeq
which tells $h$ is positive a.e..
The non-positive second variation of functional associated to inequality (\ref{shls}) gives
\beq\label{sv} \iint_{S\times S} \frac{\overline{\varphi(\zeta)}\varphi(\eta)}{|1-\zeta \cdot \bar{\eta}|^{\frac{\lambda}{2}}}d\zeta d\eta\int_S h^p-(p-1)\iint_{S\times S} \frac{h(\zeta)h(\eta)}{|1-\zeta\cdot\bar{\eta}|^{\frac{\lambda}{2}}}d\zeta d\eta\int_S h^{p-2}|\varphi|^2\le 0\eeq
for all $\varphi$ satisfying $\int_S h^{p-1}\varphi=0$.\\
\emph{Step 2:}\\
We assimilate the huge conformal symmetry group by requiring the extremizer to satisfy \emph{zero center-mass} condition. Then the opposite direction holds for the purported second-varation inequality (\ref{sv-}) derived from (\ref{sv}) by taking suitable test functions and reaches equality only by constant function, which we will prove in step 4 by computing the eigenvalues of the quardratic forms. This zero center-mass technique was introduced in \cite{hersch1970quatre} and then used by \cite{chang1987prescribing,branson2013moser,frank2012sharp,clz}. So, despite the non-uniqueness of extremizer made by the symmetry group, we can obtain a special extremizer (constant function) by repulling back any extremzier through a delicately chosen conformal transformation $\gamma$. We state this accurately in the following lemma.
\blem\label{l-zero}
For any positive extremizer of HLS inequality \emph{(\ref{shls})}, there exists a comformal transformation $\gamma: S\rightarrow S$,  s.t.
by $h\mapsto \tilde{h}=|J_{\gamma^{-1}}|^{\frac{1}{p}}h\circ \gamma^{-1}$, we get another positive extremizer $\tilde{h}$ satisfying the following zero center-mass condition:
\beq\label{zero} \int_S \zeta \tilde{h}^p(\zeta)d\zeta =0. \eeq
\elem
The proof of Lemma \ref{l-zero} is the same as that in \cite{frank2012sharp,clz}. We take a conformal map on $S\setminus \{A_\xi^{-1}o\}$:
\beq\label{gamma} \gamma^\delta_\xi=A_\xi^{-1}\circ\ca{C}\circ\ca{S}_\delta\circ\ca{C}^{-1}\circ A_\xi\eeq
with any rotation $A_\xi\in Spin(9)$ satisfying $ A_\xi^{-1}(0,1)=\xi$ (for $\xi\in S$), $\ca{S}_\delta (\delta>0)$ the dilation on $G$. Note the definition of $\gamma_\xi^\delta$ is only dependent of $\xi$, but independent of the specific choice of $A_\xi$. For any positive extremizer $h$ normalized by $\int_S h^p=1$, we take $\bb{O}^2-$valued function \[F(r\xi)=\int_S \gamma_\xi^{1-r}h^p, \quad0<r<1, \xi\in S.\] Then $F$ can be extended to a continous function on closed unit ball, with $F(\xi)=\xi$ on the boundary. From the famous Brouwer's fixed point theorem, there exists one pair $(r_0,\xi_0) (0<r_0<1, \xi_0\in S)$, s.t. $F(r_0\xi_0)=\int_S \gamma_{\xi_0}^{1-r_0}h^p=0$. Take $\gamma=\gamma_{\xi_0}^{1-r_0}$ and $\tilde{h}$ in the lemma, then $\int_S \zeta \tilde{h}^p(\zeta)d\zeta=0$ by changing variables. The lemma is proved.\\
\emph{Step 3:}\\
From Lemma \ref{l-zero}, we can assume that any positive extremizer $h$ for HLS inequality (\ref{shls}) satisfies zero center-mass condition (\ref{zero}). Then we will try to prove $h$ can only be constant function. We denote $\zeta^i\in\bb{C} (1\le i\le 4)$ to be the four complex part of $\zeta\in\bb{O}$, i.e. $\zeta=(\zeta^1,\zeta^2,\zeta^3,\zeta^4)$ with the obvious meaning.
Substituting test function $\varphi(\zeta)=h(\zeta)\zeta^i_j (1\le i\le 4, 1\le j\le 2)$ satisfying $\int_S h^{p-1}\varphi=0$,
into the second variation inequality (\ref{sv}) and summing the results by $i,j$,
we get
\[ \iint_{S\times S}\frac{h(\zeta)(\bar{\zeta}\cdot_\bb{C}\eta)h(\eta)}{|1-\zeta\cdot\bar{\eta}|^{\frac{\lambda}{2}}}d\zeta d\eta\le
(p-1)\iint_{S\times S} \frac{h(\zeta)h(\eta)}{|1-\zeta\cdot\bar{\eta}|^{\frac{\lambda}{2}}}d\zeta d\eta
\]
with $\bar{\zeta}\cdot_\bb{C}\eta=\sum_{1\le i\le 4,1\le j\le 2}\overline{\zeta^i_j}\eta^i_j$ being the complex product.
From the symmetry of left integrand on $(\zeta,\eta)$ and the basic identity
$\bar{\zeta}\cdot_\bb{C}\eta+\bar{\eta}\cdot_\bb{C}\zeta=\bar{\zeta}\cdot\eta+\bar{\eta}\cdot\zeta$,
we have
\beq\label{sv-}
\iint_{S\times S}\frac{h(\zeta)(\bar{\zeta}\cdot\eta+\bar{\eta}\cdot\zeta)h(\eta)}{|1-\zeta\cdot\bar{\eta}|^{\frac{\lambda}{2}}}d\zeta d\eta\le 2(p-1)\iint_{S\times S} \frac{h(\zeta)h(\eta)}{|1-\zeta\cdot\bar{\eta}|^{\frac{\lambda}{2}}}d\zeta d\eta.
\eeq
\emph{Step 4:}\\
By checking the second-variation inequality (\ref{sv-}), we find the opposite direction inequality holds and obtain equality if and only if $h$ is constant function.
We use octonionic analogue of Funk-Hecke formula to estimate the quadratic forms in both sides of the inequality. We list the result in the following theorem.

\bthm\label{t-isv}
\emph{[Bilinear Inequality]}\\
Let $3\le \alpha<\frac{Q}{4}, Q=22$, then for any $f$ on $S$, we have
\beq
\iint_{S\times S}\frac{\overline{f(\zeta)}(\bar{\zeta}\cdot\eta+\bar{\eta}\cdot\zeta)f(\eta)}{|1-\zeta\cdot\bar{\eta}|^{2\alpha}}d\zeta \eta\ge \frac{2\alpha}{\frac{Q}{2}-\alpha}\iint_{S\times S} \frac{\overline{f(\zeta)}f(\eta)}{|1-\zeta\cdot\bar{\eta}|^{2\alpha}}d\zeta d\eta
\eeq and $\alpha\ge 3$ is sharp, i.e. for any $\alpha<3$, there exists a function $f_\alpha$ that makes the inequality invalid. Moreover, when $\alpha>3$, equality holds if and only if $f$ is constant function;
When $\alpha=3$, equality holds if and only if $f\in W_{0,0}\bigoplus_{j\ge k\ge 2}W_{j,k}$
and if $f$ is furthermore an extremizer for HLS inequality \emph{(\ref{shls})}, $f$ can only be constant function.
\ethm
From this theorem (take $\alpha=\frac{\lambda}{4}$) and (\ref{sv-}), we see that positive extremizer satisfying zero center-mass condtion can only be be constant function. Then pulling back, all extremizers are of the form $|J_\gamma|^{\frac{1}{p}}$ ($\gamma=\gamma_\xi^\delta$ in (\ref{gamma}) for some $\delta,\xi$). Simultaneously, from the second remark following our main theorems (Theorem \ref{t-hls} and \ref{t-shls}), we know that $|J_\tau|^{\frac{1}{p}} (\tau\in Aut(S))$ (or $|J_\gamma|^{\frac{1}{p}}, \gamma=\gamma_\xi^\delta$ in (\ref{gamma}) for some $\delta,\xi$ in the proof of Lemma \ref{l-zero}) are right all the extremizers for HLS inequality on sphere (\ref{shls}). Basic computation like that in \cite{jerison1988extremals,branson2013moser,clz} gives the explicit form of all extremizers (\ref{e-s}) and all extremizers for sharp HLS inequality on group (\ref{hls}) comes from relation (\ref{cor}) and are given by (\ref{e}). The sharp constants comes from Lemma \ref{l-eig}. Concerning the proof of Theorem \ref{t-isv} involves complicated computation, we first states two lemmas about eigenvalues of operator with kernel of the form $K(\zeta\cdot\bar{\eta}), K\in L^1(B(0,1))$, where we denote $B(0,1)=\set{u\in\bb{O}:|u|<1}$ and $\partial B(0,1)=\set{u\in\bb{O}:|u|=1}\simeq\bb{S}^7$ respectively for the unit ball and sphere in $\bb{O}$.
\blem\label{l-fh}
\emph{[Quaternionic Funk-Hecke Formula]}\\
Let $K$ be a function on $B(0,1)$, the unit ball in $\bb{O}$, s.t. the following integral exists, like $K\in L^1(B(0,1))$. Then integral operator with kernel $K(\zeta\cdot\bar{\eta})$ is diagonal w.r.t bispherical harmonic decomposition \emph{(\ref{dec})}, and the eigenvalue on $(j,k)$-bispherical harmonic subspace $W_{j,k}$ is given by
\begin{align*}
\lambda_{j,k}(K)=&\frac{15\pi^4 k!}{(k+3)!}\int_0^{\frac{\pi}{2}}\cos^{j-k+7}\theta\sin^7\theta P_k^{(3,3+j-k)}(\cos2\theta)d\theta \int_{\partial B(0,1)} K(\overline{u}\cos\theta)\\
&\times\Big(a_{j,k}^0\cos(j-k)\phi+a_{j,k}^1\cos(j-k+2)\phi+a_{j,k}^2\cos(j-k+4)\phi\\
&+a_{j,k}^3\cos(j-k+6)\phi\Big)du
\end{align*} with $\re u=\cos\phi ~(\phi\in[0,\pi])$, $du$ the standard Lebesgue surface measure on $\partial B(0,1)$, the unit sphere in $\bb{O}$, $P_k^{(3,3+j-k)}(z)$ the Jacobi polynomial of order $k$ associated to weight $(1-z)^3(1+z)^{3+j-k}$  and
\begin{align*}
a_{j,k}^0=&+\frac{1}{8}\frac{1}{j-k+3}-\frac{1}{4}\frac{1}{j-k+2}+\frac{1}{8}\frac{1}{j-k+1}\\
a_{j,k}^1=&+\frac{3}{8}\frac{1}{j-k+3}-\frac{1}{4}\frac{1}{j-k+4}-\frac{1}{8}\frac{1}{j-k+1}\\
a_{j,k}^2=&-\frac{3}{8}\frac{1}{j-k+3}+\frac{1}{4}\frac{1}{j-k+2}+\frac{1}{8}\frac{1}{j-k+5}\\
a_{j,k}^3=&-\frac{1}{8}\frac{1}{j-k+3}+\frac{1}{4}\frac{1}{j-k+4}-\frac{1}{8}\frac{1}{j-k+5}.
\end{align*}
\elem
\bpf
From Schur's lemma and the irreducibility of $(j,k)$-bispherical harmonic subspace $W_{j,k}$, we see the integral operator associated to
$K(\zeta\cdot\bar{\eta})$ is diagonal with eigenvalues denoted by $\lambda_{j,k}$. Now, we compute the eigenvalues $\lambda_{j,k}$.
Assume $\{Y_{j,k}^\mu\}_{1\le\mu\le m_{j,k}}$ is a normalized orthogonal basis of $W_{j,k}$, then
\[\int_S K(\zeta\cdot\overline{\eta})Y_{j,k}(\eta)d\eta=\lambda_{j,k} Y_{j,k}(\zeta),\] and
in abuse of notation, the reproducing kernel of projection onto $W_{j,k}$ is given by
\[Z_{j,k}(\zeta,\eta)=Z_{j,k}(\zeta\cdot\bar{\eta})=\sum_{\mu=1}^{m_{j,k}}Y_{j,k}^\mu(\zeta) \overline{Y_{j,k}^\mu(\eta)},\] which coincides with (\ref{zh}) when $\eta=(0,1)$, i.e. (\ref{zh}) is invariant under the action of the subgroup fixing noth pole ($\backsimeq Spin(7)$).
Then we have
\[\int_S K(\zeta\cdot\bar{\eta})Z_{j,k}(\eta\cdot\bar{\zeta})d\eta=\lambda_{j,k}Z_{j,k}(1),\] which implies
\begin{align}\label{1}
\lambda_{j,k} &= Z^{-1}_{j,k}(1)\int_S K(\zeta\cdot\bar{\eta})Z_{j,k}(\eta\cdot\bar{\zeta})d\eta\nonumber\\
&=Z^{-1}_{j,k}(1)\int_S K(\overline{\eta_2})Z_{j,k}(\eta_2)d\eta.
\end{align}
In polar coordinates,
\begin{eqnarray*}
  \eta   &=& (\eta_1,\eta_2)\\
  \eta_1 &=& u_1\sin\theta \\
  \eta_2 &=& u_2\cos\theta\\
    u_j^8 &=& \sin\phi_j^7\sin\phi_j^6\ldots\sin\phi_j^2\sin\phi_j^1 \\
    u_j^7 &=& \sin\phi_j^7\sin\phi_j^6\ldots\sin\phi_j^2\cos\phi_j^1 \\
   \ldots & & \ldots\\
    u_j^2 &=& \sin\phi_j^7\cos\phi_j^6 \\
    u_j^1 &=& \cos\phi_j^7
\end{eqnarray*} with $\theta\in [0,\frac{\pi}{2}], \phi_j^7,\ldots, \phi_j^2\in[0,\pi],\phi_j^1\in[0,2\pi] (1\le j\le 2)$,
 we have the invariance measure
\beq\label{2} d\eta= 
 \sin^7\theta\cos^7\theta d\theta du_1du_2,
\eeq with $du_j=\sin^6\phi_{j}^7\ldots\sin\phi_{j}^2d\phi_{j}^7\ldots d\phi_{j}^2d\phi_{j}^1$.
From several formulas about special functions (hypergeometric function,  Gegenbauer and Jacobi polynomials), their relations and integral representations (15.3.4,15.3.22,15.4.5,22.10.11,22.5.42 in \cite{abramowitz2012handbook}), we write the explicit formula of zonal harmonics (\ref{zh}) to the following form (here we ignore some dispensable constant, see (\ref{1})),
\begin{align}\label{13}
  Z_{j,k}(\eta_2)\sim &\sin^{-5}\phi\bigg[\frac{\sin(j-k+5)\phi+\sin(j-k+1)\phi}{4(j-k+3)}+\frac{\sin(j-k+3)\phi}{j-k+3}\nonumber\\
& -\frac{1}{2}\bra{\frac{\sin(j-k+3)\phi+\sin(j-k+1)\phi}{j-k+2}+\frac{\sin(j-k+5)\phi+\sin(j-k+3)\phi}{j-k+4}}\nonumber\\
& +\frac{1}{4}\bra{\frac{\sin(j-k+1)\phi}{j-k+1}+\frac{\sin(j-k+5)\phi}{j-k+5}}
\bigg] \cos^{j-k}\theta~ P_k^{(3,3+j-k)}(\cos2\theta),
\end{align} as we can simplify the two terms in (\ref{zh}) into
\begin{align*}
&\cos^{j-k}\phi~{}_2F_1\bra{-\frac{j-k}{2},\frac{1-(j-k)}{2};\frac{7}{2};-\tan^2\phi}\\
&={}_2F_1\bra{-\frac{j-k}{2},3+\frac{j-k}{2};\frac{7}{2};\sin^2\phi}\\
&={}_2F_1\bra{-(j-k),6+(j-k);\frac{7}{2};\sin^2\frac{\phi}{2}}\\
&=\frac{5!(j-k)!}{(j-k+5)!} C_{j-k}^{(3)}(\cos\phi)\\
&=\frac{15}{2}\sin^{-5}\phi\int_{0}^{\phi}\cos(j-k+3)t~(\cos t-\cos\phi)^2dt\\
&=\frac{15}{2}\sin^{-5}\phi\bigg[\frac{\sin(j-k+5)\phi+\sin(j-k+1)\phi}{4(j-k+3)}+\frac{\sin(j-k+3)\phi}{j-k+3}\\
&\quad~ -\frac{1}{2}\bra{\frac{\sin(j-k+3)\phi+\sin(j-k+1)\phi}{j-k+2}+\frac{\sin(j-k+5)\phi+\sin(j-k+3)\phi}{j-k+4}}\\
&\quad~ +\frac{1}{4}\bra{\frac{\sin(j-k+1)\phi}{j-k+1}+\frac{\sin(j-k+5)\phi}{j-k+5}}
\bigg],\\
&\cos^{j+k}\theta~{}_2F_1\bra{-k,\frac{-j-6}{2};4;-\tan^2\theta}\\
&=\cos^{j-k}\theta~{}_2F_1\bra{-k,7+j;4;\sin^2\theta}\\
&=\frac{3!k!}{(k+3)!} \cos^{j-k}\theta~ P_k^{(3,3+j-k)}(\cos2\theta),
\end{align*} where $C_{j-k}^{(3)},  P_k^{(3,3+j-k)}$ are respectively Gegenbauer and Jacobi polynomials.
Putting (\ref{13}) and (\ref{2}) into (\ref{1}) and noting that $\frac{5!(j-k)!}{(j-k+5)!} C_{j-k}^{(3)}(1)=1,  \frac{3!k!}{(k+3)!}P_k^{(3,3+j-k)}(1)=1$ (22.2.1 and 22.2.3 in \cite{abramowitz2012handbook}), we have
\begin{align*}
\lambda_{j,k}(K)=&\frac{15\pi^4 k!}{(k+3)!}\int_0^{\frac{\pi}{2}}\cos^{j-k+7}\theta\sin^7\theta P_k^{(3,3+j-k)}(\cos2\theta)d\theta \int_{\partial B(0,1)} K(\overline{u_2}\cos\theta)\\
&\times\Big(a_{j,k}^0\cos(j-k)\phi_2^7+a_{j,k}^1\cos(j-k+2)\phi_2^7+a_{j,k}^2\cos(j-k+4)\phi_2^7\\
&+a_{j,k}^3\cos(j-k+6)\phi_2^7\Big)du_2
\end{align*} with $a_{j,k}^i (0\le i\le 3)$ defined in the lemma, then lemma is proved.
 \epf
In order to prove Theorem \ref{t-isv}, it suffices to compute eigenvalues of functions of two forms
$K^\alpha_1(z)=|1-z|^{-2\alpha},K^\alpha_2(z)=|z|^2|1-z|^{-2\alpha}$, noting that
$\zeta\cdot\bar{\eta}+\eta\cdot\bar{\zeta}=2\re\zeta\cdot\bar{\eta}=1+|\zeta\cdot\bar{\eta}|^2-|1-\zeta\cdot\bar{\eta}|^2$.
\blem\label{l-eig}\emph{[Eigenvalues]}\\
Given $-1<\alpha<\frac{Q}{4}$, denote \emph{(15.1.1 and 15.1.20 in \cite{abramowitz2012handbook})}
\beq\label{A}
A(a,b,c) \triangleq \sum_{\mu\ge 0}\frac{\Gamma(\mu+a)\Gamma(\mu+b)}{\mu!\Gamma(\mu+c)}=\frac{\Gamma(a)\Gamma(b)\Gamma(c-a-b)}{\Gamma(c-a)\Gamma(c-b)}
\eeq for any $c>a+b$ and \beq\label{abc}(a,b,c)= (j+\alpha,k+\alpha-3,j+k+\frac{Q}{2}-3)\eeq particularly in the following part of this paper.
Then\\
\emph{(1)} The eigenvalues of integral operators associated to kernel $K^\alpha_1(\zeta\cdot\bar{\eta})=|1-\zeta\cdot\bar{\eta}|^{-2\alpha}$ are given by
\begin{align}
\lambda_{j,k}(K^\alpha_1)&=~\frac{2\pi^8}{\Gamma(\alpha)\Gamma(\alpha-3)} A(a,b,c)\nonumber\\
&=~\frac{2\pi^8\Gamma(\frac{Q}{2}-2\alpha)}{\Gamma(\alpha)\Gamma(\alpha-3)}
\frac{\Gamma(j+\alpha)}{\Gamma(j+\frac{Q}{2}-\alpha)}
\frac{\Gamma(k+\alpha-3)}{\Gamma(k+\frac{Q}{2}-\alpha-3)}.
\end{align}
\emph{(2)} The eigenvalues of integral operators associated to kernel $K^\alpha_2(\zeta\cdot\bar{\eta})=|\zeta\cdot\bar{\eta}|^2|1-\zeta\cdot\bar{\eta}|^{-2\alpha}$ are given by
\beq\lambda_{j,k}(K^\alpha_2)=C_{j,k}^{\alpha}\lambda_{j,k}(K^\alpha_1)\eeq with
\begin{align}
C_{j,k}^{\alpha}=&1-(\alpha-4)(c-a-b)\nonumber\\
&\times\left(\frac{1}{(a-1)(c-a)}+\frac{1}{(b-1)(c-b)}-(\alpha-4)\frac{c-a-b+1}{(a-1)(b-1)(c-a)(c-b)}\right)\nonumber\\
=&1-\frac{(\alpha-4)(\frac{Q}{2}-2\alpha)}{(j+\alpha-1)(k+\alpha-4)(k+\frac{Q}{2}-\alpha-3)(j+\frac{Q}{2}-\alpha)} \nonumber\\
&\times\Big(-(j+\alpha)^2-(k+\alpha-3)^2+(j+k+\frac{Q}{2}-2)(j+2\alpha+k-3)\nonumber\\
&-2(j+k+\frac{Q}{2}-3)-(\alpha-4)(\frac{Q}{2}-2\alpha+1)\Big).
\end{align}  In the singular point $\alpha=0,1,2,3,4$, the above formula can be viewed as limit, fixing $j,k$.
\elem
\bpf
Taking $K=K_1^\alpha, K_2^\alpha$ in last Lemma \ref{l-fh}, we get
\begin{align}
\lambda_{j,k}(K_1^\alpha) = &\frac{16\pi^7 k!}{(k+3)!}\int_0^{\frac{\pi}{2}} \cos^{j-k+7}\theta\sin^7\theta P_k^{(3,3+j-k)}(\cos2\theta)d\theta\nonumber\\
&\times \int_{0}^{\pi} (1-2\cos\theta\cos\phi+\cos^2\theta)^{-\alpha} \Big(a_{j,k}^0\cos(j-k)\phi+a_{j,k}^1\cos(j-k+2)\phi\nonumber\\
&+a_{j,k}^2\cos(j-k+4)\phi+a_{j,k}^3\cos(j-k+6)\phi\Big)d\phi, \label{3-1}\\
\lambda_{j,k}(K_2^\alpha) = &\frac{16\pi^7 k!}{(k+3)!}\int_0^{\frac{\pi}{2}} \cos^{j-k+9}\theta\sin^7\theta P_k^{(3,3+j-k)}(\cos2\theta)d\theta\nonumber\\
&\times \int_{0}^{\pi} (1-2\cos\theta\cos\phi+\cos^2\theta)^{-\alpha} \Big(a_{j,k}^0\cos(j-k)\phi+a_{j,k}^1\cos(j-k+2)\phi\nonumber\\
&+a_{j,k}^2\cos(j-k+4)\phi+a_{j,k}^3\cos(j-k+6)\phi\Big)d\phi. \label{3-2}
\end{align}
Using Gegenbauer polynomials, the following fact exists ((5.11) in \cite{frank2012sharp})
\begin{align}\label{4}
&\int_0^\pi d\phi (1+\cos^2\theta-2\cos\phi\cos\theta)^{-\alpha}\cos(j-k)\phi\nonumber\\
&=\frac{\pi}{\Gamma^2(\alpha)}\sum_{\mu\ge 0}\cos^{|j-k|+2\mu}\theta\frac{\Gamma(\mu+\alpha)\Gamma(\mu+|j-k|+\alpha)}{\mu!(\mu+|j-k|)!}.
\end{align}
Putting (\ref{4}) into (\ref{3-1},\ref{3-2}), we get
\begin{align}
\lambda_{j,k}(K_1^\alpha)=& \frac{16\pi^8 k!}{(k+3)!\Gamma^2(\alpha)}\sum_{i=0}^3 a_{j,k}^i
\sum_{\mu\ge 0}\frac{\Gamma(\mu+\alpha)\Gamma(\mu+|j-k|+2i+\alpha)}{\mu!(\mu+|j-k|+2i)!}\nonumber\\
&\times \int \cos^{2|j-k|+7+2(i+\mu)}\theta\sin^7\theta P_k^{(3,3+j-k)}(\cos2\theta) d\theta, \label{5-1}\\
\lambda_{j,k}(K_2^\alpha)=& \frac{16\pi^8 k!}{(k+3)!\Gamma^2(\alpha)}\sum_{i=0}^3 a_{j,k}^i
\sum_{\mu\ge 0}\frac{\Gamma(\mu+\alpha)\Gamma(\mu+|j-k|+2i+\alpha)}{\mu!(\mu+|j-k|+2i)!}\nonumber\\
&\times \int \cos^{2|j-k|+7+2(i+\mu+1)}\theta\sin^7\theta P_k^{(3,3+j-k)}(\cos2\theta) d\theta. \label{5-2}
\end{align}
Using Rodrigues' formula (22.2.1 in \cite{abramowitz2012handbook}) \[P_k^{(3,j-k+3)}(t)=\frac{(-1)^k}{2^kk!}(1-t)^{-3}(1+t)^{-(j-k+3)}\frac{d^k}{dt^k}\{(1-t)^{k+3}(1+t)^{j+3}\},\]
changing variable $\cos2\theta=t$ and integrating by part, we get
\begin{align}\label{6}
 &\int \cos^{2(j-k)+7+2\mu}\theta\sin^7\theta P_k^{(3,3+j-k)}(\cos2\theta) d\theta \nonumber\\
 &=\frac{(-1)^k}{2^{j+8+\mu}k!}\int_{-1}^1 (1+t)^\mu\frac{d^k}{dt^k}\set{(1-t)^{3+k}(1+t)^{3+j}} dt\nonumber\\
 &=\chi_{\mu\ge k}\frac{\mu!}{2^{\mu+j+8}k!(\mu-k)!}\int_{-1}^1(1+t)^{\mu+j-k+3}(1+t)^{3+k}dt\nonumber\\
 &=\chi_{\mu\ge k}\frac{\mu!B(\mu+j-k+4,k+4)}{2(\mu-k)!k!}\nonumber\\
 &=\chi_{\mu\ge k}\frac{(k+3)!}{2k!}\frac{\mu!(\mu+j-k+3)!}{(\mu-k)!(\mu+j+7)!}.
\end{align}
Assume $k\ge 3$ for $K_1^\alpha$ (and $k\ge 4$ for $K_2^\alpha$) and put (\ref{6}) into (\ref{5-1},\ref{5-2}), then we have
\begin{align}
\lambda_{j,k}(K_1)=& \frac{8\pi^8}{\Gamma^2(\alpha)}\sum_{\mu\ge k}\frac{\mu!(\mu+j-k+3)!}{(\mu-k)!(\mu+j+7)!}
\sum_{i=0}^3 a_{j,k}^i \frac{\Gamma(\mu-i+\alpha)\Gamma(\mu+|j-k|+i+\alpha)}{(\mu-i)!(\mu+|j-k|+i)!}, \label{7-1}\\
\lambda_{j,k}(K_2)=& \frac{8\pi^8}{\Gamma^2(\alpha)}\sum_{\mu\ge k}\frac{\mu!(\mu+j-k+3)!}{(\mu-k)!(\mu+j+7)!}
\sum_{i=0}^3 a_{j,k}^i \frac{\Gamma(\mu-1-i+\alpha)\Gamma(\mu-1+|j-k|+i+\alpha)}{(\mu-1-i)!(\mu-1+|j-k|+i)!} \label{7-2}.
\end{align}
Noting that
\begin{align*}
a_{j,k}^0=&\frac{1}{8}\bra{\frac{1}{j-k+3}+\frac{-2}{j-k+2}+\frac{1}{j-k+1}}= +\frac{1}{4}\frac{1}{(j-k+3)(j-k+2)(j-k+1)},\\
a_{j,k}^1=&\frac{1}{8}\bra{\frac{3}{j-k+3}+\frac{-2}{j-k+4}+\frac{-1}{j-k+1}}=-\frac{3}{4}\frac{1}{(j-k+3)(j-k+4)(j-k+1)},\\
a_{j,k}^2=&\frac{1}{8}\bra{\frac{-3}{j-k+3}+\frac{2}{j-k+2}+\frac{1}{j-k+5}}=+\frac{3}{4}\frac{1}{(j-k+3)(j-k+2)(j-k+5)},\\
a_{j,k}^3=&\frac{1}{8}\bra{\frac{-1}{j-k+3}+\frac{2}{j-k+4}+\frac{-1}{j-k+5}}=-\frac{1}{4}\frac{1}{(j-k+3)(j-k+4)(j-k+5)},
\end{align*}
and through easy but boring computation we have
\begin{align}\label{8}
&\sum_{i=0}^3 a_{j,k}^i \frac{\Gamma(\mu-i+\alpha)\Gamma(\mu+|j-k|+i+\alpha)}{(\mu-i)!(\mu+|j-k|+i)!}\nonumber\\
&=\frac{1}{4} (\alpha-1)(\alpha-2)(\alpha-3)\frac{\Gamma(\mu+j-k+\alpha)\Gamma(\mu+\alpha-3)}{\mu!(\mu+j-k+3)!}.
\end{align}
So putting (\ref{8}) into (\ref{7-1},\ref{7-2}) and using (\ref{A},\ref{abc}) we get
\begin{align}
\lambda_{j,k}(K_1)=& \frac{2\pi^8}{\Gamma(\alpha)\Gamma(\alpha-3)} \sum_{\mu\ge k}\frac{\mu!(\mu+j-k+3)!}{(\mu-k)!(\mu+j+7)!} \frac{\Gamma(\mu+j-k+\alpha)\Gamma(\mu+\alpha-3)}{\mu!(\mu+j-k+3)!}\nonumber\\
=&  \frac{2\pi^8}{\Gamma(\alpha)\Gamma(\alpha-3)} \sum_{\mu\ge k} \frac{\Gamma(\mu+j-k+\alpha)\Gamma(\mu+\alpha-3)}{(\mu-k)!(\mu+j+7)!}\nonumber\\
=& \frac{2\pi^8\Gamma(\frac{Q}{2}-2\alpha)}{\Gamma(\alpha)\Gamma(\alpha-3)} A(a,b,c)\nonumber\\
=& \frac{2\pi^8\Gamma(\frac{Q}{2}-2\alpha)}{\Gamma(\alpha)\Gamma(\alpha-3)}
\frac{\Gamma(j+\alpha)}{\Gamma(j+\frac{Q}{2}-\alpha)}
\frac{\Gamma(k+\alpha-3)}{\Gamma(k+\frac{Q}{2}-\alpha-3)}, \label{9-1}\\
\lambda_{j,k}(K_2)=&  \frac{2\pi^8}{\Gamma(\alpha)\Gamma(\alpha-3)} \sum_{\mu\ge k}\frac{\mu!(\mu+j-k+3)!}{(\mu-k)!(\mu+j+7)!} \frac{\Gamma(\mu-1+j-k+\alpha)\Gamma(\mu-1+\alpha-3)}{(\mu-1)!(\mu-1+j-k+3)!}\nonumber\\
=&  \frac{2\pi^8}{\Gamma(\alpha)\Gamma(\alpha-3)} \sum_{\mu\ge 0} \frac{\Gamma(\mu+j+\alpha)\Gamma(\mu+k+\alpha-3)}{\mu!(\mu+j+k+7)!}
\frac{(\mu+k)(\mu+j+3)}{(\mu+j+\alpha-1)(\mu+k+\alpha-4)}. \label{9-2}
\end{align}
Noting that
\begin{align*}
&\frac{(\mu+k)(\mu+j+3)}{(\mu+j+\alpha-1)(\mu+k+\alpha-4)}=~\frac{(\mu+a-1-(\alpha-4))(\mu+b-1-(\alpha-4))}{(\mu+a-1)(\mu+b-1)}\nonumber\\
&= 1-(\alpha-4)\left(\frac{1}{\mu+a-1}+\frac{1}{\mu+b-1}-(\alpha-4)\frac{1}{(\mu+a-1)(\mu+b-1)}\right)
\end{align*} and using (\ref{A}) again
we have
\beq\label{10}\lambda_{j,k}(K^\alpha_2)=C_{j,k}^{\alpha}\lambda_{j,k}(K^\alpha_1)\eeq with
\begin{align*} C_{j,k}^{\alpha}=& 1-(\alpha-4)(c-a-b)\\
& \times \left(\frac{1}{(a-1)(c-a)}+\frac{1}{(b-1)(c-b)}-(\alpha-4)\frac{c-a-b+1}{(a-1)(b-1)(c-a)(c-b)}\right).
\end{align*}
For $k<3$ ($k<4$), it's easy to check that above formulas (\ref{9-1}) ((\ref{9-2})and (\ref{10})) still holds.
\epf
We now proceed the \emph{proof of Theorem} \ref{t-isv} using Lemma \ref{l-eig}:\\
\bpf
Let $f=\sum_{j\ge k\ge 0}f_{j,k} \,(f_{j,k}\in W_{j,k})$ be the bispherical harmonic decomposition of $f$ w.r.t (\ref{dec}). From $\zeta\cdot\bar{\eta}+\eta\cdot\bar{\zeta}=1+|\zeta\cdot\bar{\eta}|^2-|1-\zeta\cdot\bar{\eta}|^2$, to prove the bilinear inequality in Theorem \ref{t-isv}, it suffices to prove
\beq\label{12}
\sum_{j\ge k\ge 0} \bra{\lambda_{j,k}(K_1^\alpha)+\lambda_{j,k}(K_2^\alpha)-\lambda_{j,k}(K_1^{\alpha-1})}\norm{f_{j,k}}_2^2 \ge \sum_{j\ge k\ge 0}\frac{2\alpha}{\frac{Q}{2}-\alpha}\lambda_{j,k}(K_1^\alpha)\norm{f_{j,k}}_2^2.
\eeq
From Lemma \ref{l-eig}, we have $\lambda_{j,k}(K_1^\alpha)>0$ when $\alpha>3$ and when $\alpha=3, \lambda_{j,k}(K_1^\alpha)=0$ unless $k=0$ and we also have relation
\[\lambda_{j,k}(K^{\alpha-1}_1)=\lambda_{j,k}(K^{\alpha}_1)(\alpha-1)(\alpha-4)\frac{(c-a-b)(c-a-b+1)}{(a-1)(b-1)(c-a)(c-b)}.\]
Then (\ref{12}) holds if we can prove the following inequality:\\
When $\alpha>3$,
\begin{align*}
&2-(\alpha-4)(c-a-b)\left(\frac{1}{(a-1)(c-a)}+\frac{1}{(b-1)(c-b)}+\frac{3(c-a-b+1)}{(a-1)(b-1)(c-a)(c-b)}\right)\\
&\ge \frac{2\alpha}{\frac{Q}{2}-\alpha},
\end{align*} which is
\[(\alpha-4)\left(\frac{1}{(a-1)(c-a)}+\frac{1}{(b-1)(c-b)}+\frac{3(c-a-b+1)}{(a-1)(b-1)(c-a)(c-b)}\right)\le \frac{2}{\frac{Q}{2}-\alpha}\] since $c-a-b=\frac{Q}{2}-2\alpha>0$.\\
For $\alpha-4> 0$, then it suffices to check
\[\frac{2}{(\alpha-4)(\frac{Q}{2}-\alpha)}-
\left(\frac{1}{(a-1)(c-a)}+\frac{1}{(b-1)(c-b)}\right)\ge \frac{3(c-a-b+1)}{(a-1)(b-1)(c-a)(c-b)}\]
Substituting $(a,b,c)$ in (\ref{abc}), the inequality becomes
\begin{align}\label{11}
&\frac{j(\alpha-4)+k(\frac{Q}{2}-\alpha)+kj}{(\alpha-4)(\frac{Q}{2}-\alpha)(\alpha-4+k)(\frac{Q}{2}-\alpha+j)}
+\frac{(k-3)(\alpha-4)+(j+3)(\frac{Q}{2}-\alpha)+(j+3)(k-3)}{(\alpha-4)(\frac{Q}{2}-\alpha)(\alpha-4+j+3)(\frac{Q}{2}-\alpha+k-3)}\nonumber\\
&\ge \frac{3(\frac{Q}{2}-2\alpha+1)}{(\alpha-4+k)(\frac{Q}{2}-\alpha+j)(\alpha-4+j+3)(\frac{Q}{2}-\alpha+k-3)}.
\end{align}
It's easy to check this noting that the denominators in left side are less than that in right side, while the sum of two numerators in left side is bigger than that in the right side. Moreover, the inequality reaches equality if and only if $j=k=0$.
For $3<\alpha< 4$, the opposite of (\ref{11}) holds by checking signs of the numerators and denominators and noting that it still holds that the sum of numerators in left side is bigger than that in the right side. Inequality again reaches equality only when $j=k=0$. $\alpha=4$ case can be checked similarly or implied from limitation argument. Note that, when $\alpha=3$, for $k=0$ the opposite of (\ref{11}) holds and reaches equality if and only if $j=0$, for $k>1$, $\frac{\Gamma(k+\alpha-3)}{\Gamma(\alpha-3)}$ multiple of opposite of (\ref{11}) reaches equality as every eigenvalue vanishes, for $k=1$, it's easy to check that $\frac{\Gamma(k+\alpha-3)}{\Gamma(\alpha-3)}$ multiple of opposite of (\ref{11}) holds strictly. In conclusion, when $\alpha=3$, (\ref{12}) holds and reaches equality if and only if $f$ is in direct sum of subspaces $W_{j,k}$ with $j=k=0$ or $k>1$, then from Euler-Lagrange equation or original functional, we see extremizer satisfying zero center-mass condition still can only be constant function (in $W_{0,0}$), see \cite{clz}. For $\alpha<3$, (\ref{12}) and therefore the inequality in Theorem \ref{t-isv},  fails in some subspaces, which tells our method doesn't work.
\epf

\section*{Acknowledgement}
The work was partially done during An Zhang's visit as a joint Ph.D. student in the department of mathematics, UC, Berkeley. He would like to thank especially the department for hospitality.

\bibliographystyle{amsalpha}
\bibliography{hls}

\vskip 3\baselineskip
\flushleft

Michael Christ \\
Department of Mathematics,
University of California, Berkeley;
Berkeley, CA, USA.\\
\emph{email:} \texttt{mchrist@math.berkeley.edu}
\vskip 2\baselineskip
Heping Liu\\
School of Mathematical Science,
Peking University;
Beijing, China.\\
\emph{email:} \texttt{hpliu@math.pku.edu.cn}
\vskip 2\baselineskip
An Zhang\\
School of Mathematical Science,
Peking University;
Beijing, China.\\
\emph{email:} \texttt{anzhang@pku.edu.cn}\\
Current Address:\\
Department of Mathematics,
University of California, Berkeley;
Berkeley, CA, USA. \\
\emph{email:} \texttt{anzhang@math.berkeley.edu}

\end{document}